# シェアサイクルにおける自転車の再配置作業ルート最適化に関する研究

Qiwei Luo
Graduate School of Information, Production and Systems, Waseda University

**要旨**：シェアサイクルのネットワークでは、需要と供給のバランスが崩れている。その結果、共有資源が無駄になり、事業者の管理コストが増大する。そこで本文では、東京都港区のシェアサイクルのポート点の容量、運搬車の実効性能、空間分布の特性から、サイクルの再配置問題を積載量制限と実際の往復距離のある運搬経路問題（注：VRP問題，固有名詞）に変換し、実際のエネルギー消費量を考慮したシェアサイクル再配置経路の最適化モデルを構築する。最後に、モデルの解を求めるために異なる貪欲アルゴリズムを使用し、アルゴリズムの再配置効果を比較し、シェアサイクルの再配置問題に対して、実データに基づく再配置の最適化されたフレームワークと改良されたアルゴリズムを提案する。
**キーワード**：シェアサイクル，再配置の最適化，最優経路

# 1 紹介

シェアサイクルとは、レンタサイクル事業者がキャンパスや地下鉄、ビジネスエリアなどで自転車をタイムシェアリングするサービスである。シェアリングエコノミーのブームの産物として、シェア自転車は都市の公共交通においてますます重要な役割を果たしており、都市部の利用者に、より環境に優しい移動手段を提供し、大気汚染も著しく減少させている。

シェアサイクルの再配置とは、いくつかのアルゴリズムを通じて、ネットワーク内の自転車資源を合理的に配分することであり、本文はシェアサイクルの再配置の方法を最適化し、研究するものである。その貢献は、(1)事業者にとっては、再配置の経路の最適化により、運搬車が業務中に移動する距離を短縮することができ、再配置の業務を各ポート点に合理的に割り当てることにより、運搬車と関連スタッフの投入を減らし、コストを削減することができる；(2)シェアサイクルを利用するユーザーにとっては、再配置の最適化により、できるだけ短時間でシェアサイ

クルを利用できるようにする。本文では、配送開始点をどのように選択するかも重要な問題であり、配送過程における運搬車の実際の状況をシミュレートするために、通常のクラスタリングアルゴリズムを用いて仮想配送点を選択するのではなく、実際のポート点を出発点として使用するのだ；(3) 社会資源に対して、シェアサイクル資源を効果的に配分することで、社会資源の利用率を向上させることができる。シェアサイクルの再配置を最適化することで、サービスエリアにおけるシェアサイクルの需要を把握し、最適な再配置経路を通じて、限られた資源でエリアの社会的需要を満たすことができる。本文は、事業者に基づくシェアサイクル再配置の最適化に焦点を当てる。その本質は運搬経路問題である。（Ｖｅｈｉｃｌｅ　Ｒｏｕｔｉｎｇ　Ｐｒｏｂｌｅｍ,ＶＲＰ），ＤａｎｔｚｉｇとＲａｍｓｅｒが１９５９年に提出された運搬経路問題（ＶＲＰ），それは、一定のユーザーを指し、それぞれ異なる量の需要を持ち、物流センターはユーザーに物品を提供し、車両隊は物品の配給を担当および適切な経路を計画する。目標はユーザーのニーズを満たすようにすることであり、最短距離、最小コスト、最小消費時間などを一定の制約の下で達成することである。シェアサイクルの再配置における"ユーザー"とはポート点である。それはサイクルの需要側であると同時に、供給側でもある。そのため、運搬車は配送を完了させながら、商品をピックアップするタスクを完了させる必要がある。

　　本文では、公開されているシェアサイクルポート点の緯度と経度の座標を用いて、東京都港区内のシェアサイクルのポート点のデータセッ

トを構築した。配車トラックの性能、シェアサイクルの重量、電気代のデータに基づいてコスト推定アルゴリズムを設計した。港区の実際のポート点を再配置の対象として、最短距離、シェアサイクルの最大需要、および目的関数としてVRPモデルをそれぞれ構築した。再配置のモデルを解くために、異なる貪欲なアルゴリズムが設計し、それぞれのアルゴリズムの流れを詳細に説明した。最後に、いくつかの現実的なルールに従い、シェアサイクルの需要データと現在の時間帯における自転車ネットワークの分布を生成するためのデータシミュレーションを実施する。データのシミュレーションによって、実際の再配置でモデルがどのように機能するかを確認する。

## 2 データの準備

### 2.1 シェアサイクルのポート点の配置

研究エリアは東京の南東に位置し、東京湾に近い港区である。NTTドコモが公開しているシェアサイクルのポート点のデータセットから、150箇所のポート点の緯度と経度を取得した。地図上のすべてのポート点を検索することで、データセットから地下駅に近いポート点11カ所と、学校や会社に隣接するポート点23カ所を特定した。このデータセットの1行目にランダムに配置された点を起点としたのが、六本木ファーストビルである。その経度は139.741476、緯度は35.66241である。その点のグラフはいかにある図1に示したように。

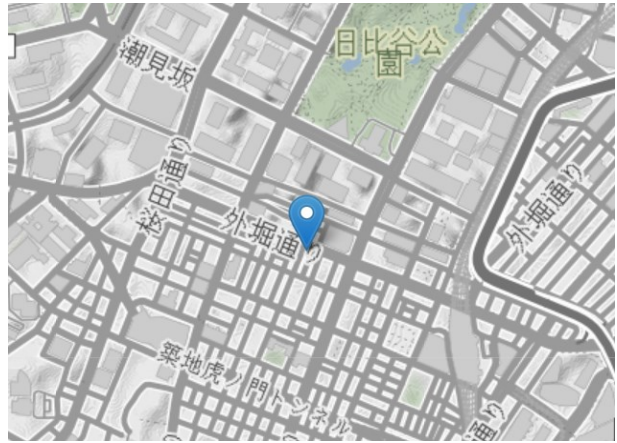

図1 六本木ファーストビル

この150カ所を、縦軸を経度、横軸を緯度としたシェアサイクルのポート点配置マップにした。次の図2のようになる。それぞれのデータポイントには、名前、経度、緯度、ラッコ数などの情報が含まれている。

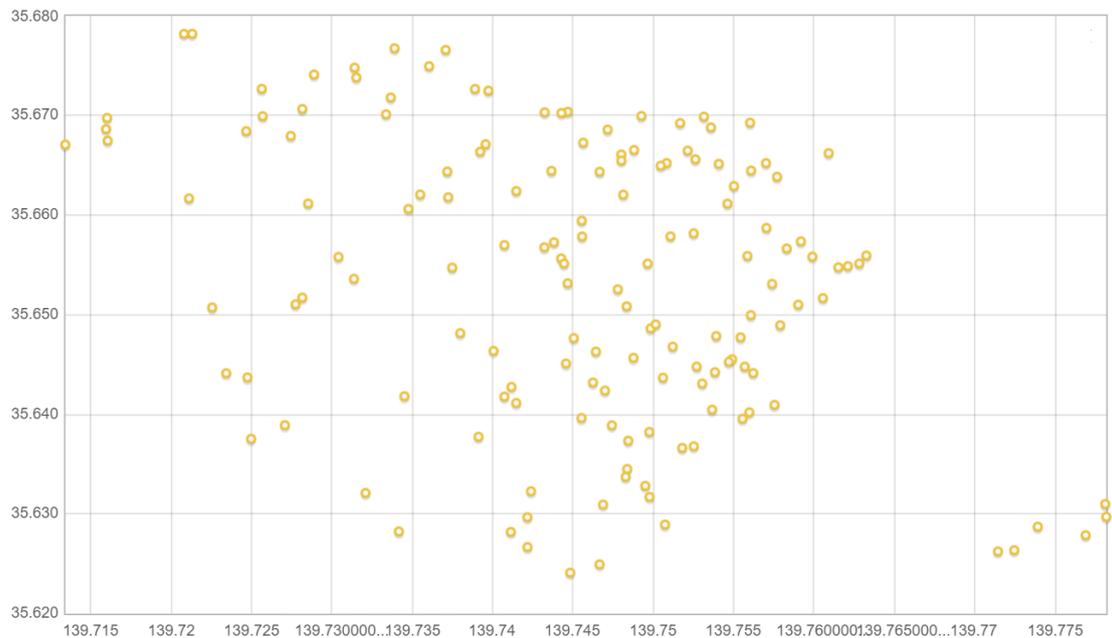

図2 シェアサイクルのポート点配置図

2.2 距離行列

Openrouteserviceは無料の地理データを提供できるウェブサイトで、Openrouteserviceの API 経由で 2 つのデータポイント間の実際の距離を調べだした、その往復距離は異なっている。Openrouteserviceでは 50 × 50 の距離行列しか解けないので、次元の異なる十数の距離行列をつなぎ合わせる必要がある、最後に 150 × 150 の行列を作り出す。距離行列の次元が高いと表示しにくいので、150 次元の距離行列をヒートマップとして表示する。横軸と縦軸はデータポイントの識別子を表す。例えば、横軸 1 と縦軸 2 はデータポイント 1 とデータポイント 2 の往復データを表す。右側の色軸は距離の大きさを表す。

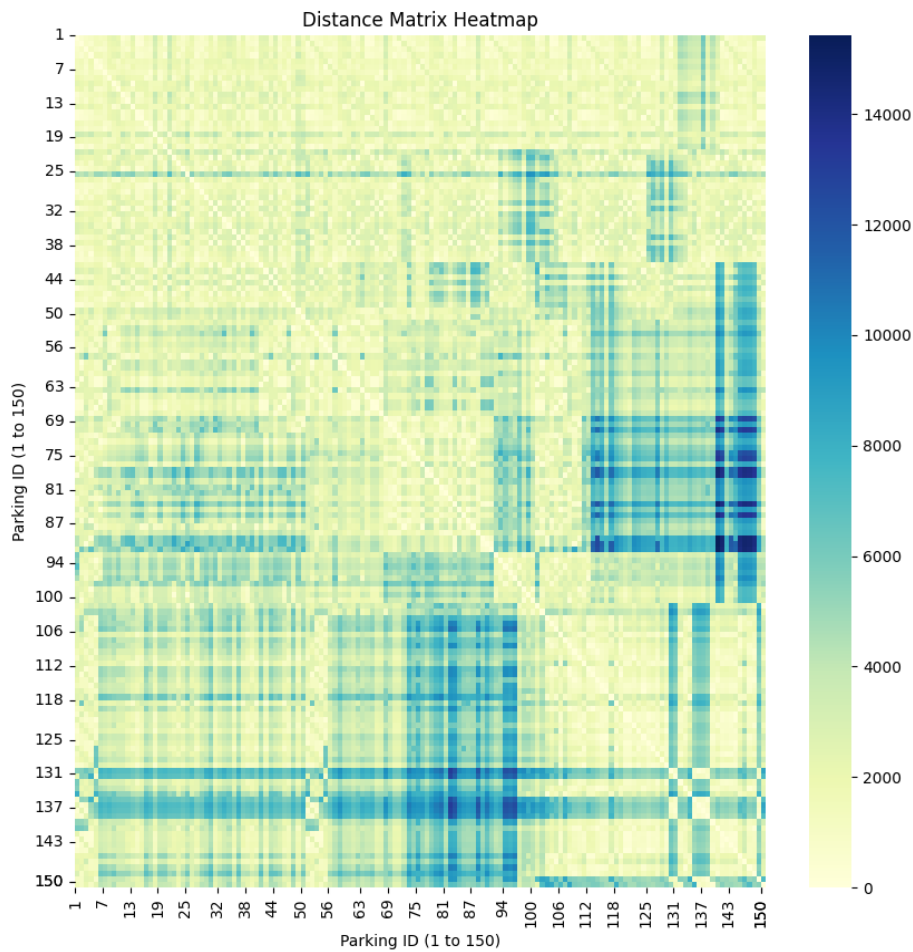

図 3 距離行列のヒートマップ

2.3 シェアサイクルと運搬車の性能指数

　　ドコモの公式サイトの情報によると、シェア自転車の車種は Bikke 2e、流通トラックの最新車種は F1 TRUCK とのこと。これら 2 台の自転車とトラックの具体的な性能パラメータは、folofly の公式ウェブサイトと Amazon で知ることができます。これらのパラメーターは、後で制約としてスケジューリング モデルに統合されます。N/A は、関連する情報がないことを意味します。表 1 に示すとおり。

表 1 シェア自転車と物流トラックの性能パラメータ

|  | シェアサイクル | 運搬車 |
| --- | --- | --- |
| 品番 | Docomo | Folofly |
| 生産会社 | Bikke 2e | F1 TRUCK |
| 重さ | 30.2 KG | 1320 KG |
| 積載量 | N/A | 1150 KG |
| 航続距離 | 38 km | 300 km |
| バッテリー容量 | 219 Wh | 38.7 kwh |

2.4 電気料金

　　GlobalPetrolPrices.com は、150 カ国以上の自動車燃料、電気、天然ガスの小売価格を掲載するウェブサイトです。サイト上の各ポート点のデータは手作業で収集されており、この分野ではリーダー的存在だ。GlobalPetrolPrices.com で確認したところ、2022 年 12 月の日本の企業向け電気料金は 1kWh あたり 0.234 ドルで、円に換算すると下の図に示

したように 33.56 円となる。この料金は、運搬車の輸送コストの計算に影響する。

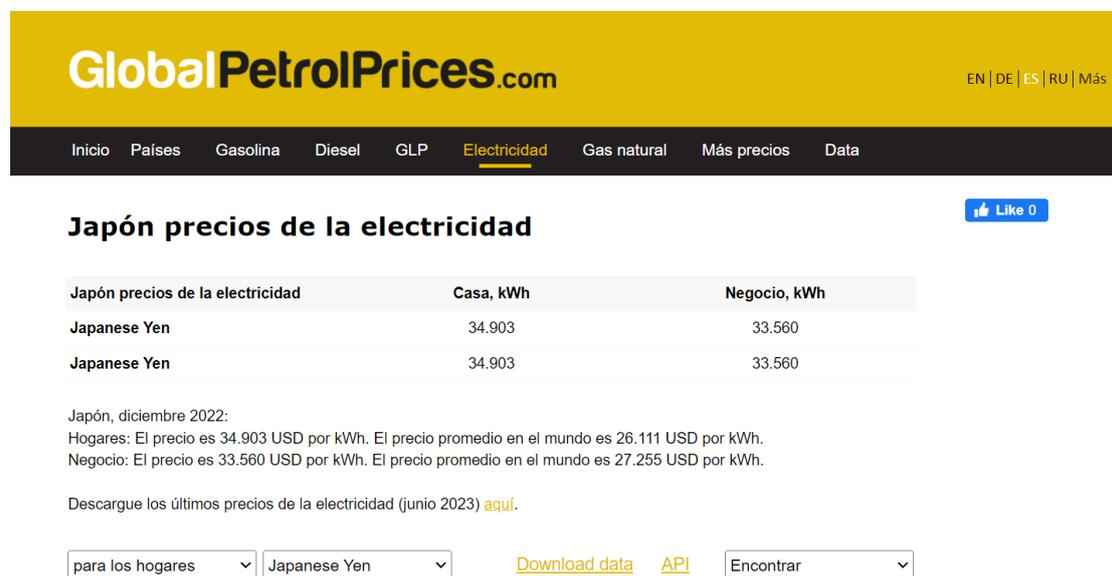

図 4 日本 2022 年 12 月の電気料金

2.5 シェアサイクルの使用率

日本自動車産業振興協会は、600 以上のサンプルを分析し、シェアサイクルの使用率と時間帯の関係についての研究データを発表している。この研究データは、その後の時間変化に応じてサンプルデータを補正する際に参考にする。次の図 5 のようになる。

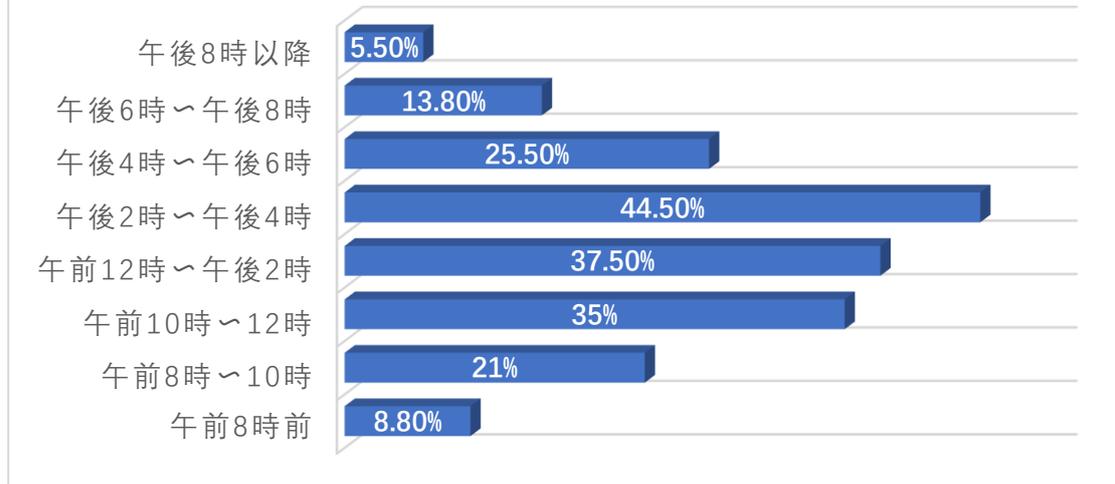

図5 日本シェアサイクルの使用率

3 再配置経路のモデル

3.1 実験の設計

　シェアサイクルの最大台数のあるパーセンテージをベースラインとし、実際のポート点の数がベースラインを超えた場合、それを「オーバーフロー」と判定する。実際のポート点の数がベースラインより少ない場合は、それを「アンダーフロー」と断定する。「オーバーフロー」も「アンダーフロー」も、ポート点に最大数のシェアサイクルが用意されるまで、運搬車によって再配分される必要がある。車両は出発点から配車エリア内のシェアサイクルのポート点を通過してサービスを提供するが、すでに需要が達成されたポート点を繰り返し訪問することはない。そして、ポート点ごとにシェアサイクルの台数がベースラインを超えると、余ったシェアサイクルを運び去る。そのポート点のシェアサイクル

の台数がベースラインより多ければ、その分だけのシェアサイクルを運搬車から降ろして補充する。運搬車の積載重量とバッテリー容量の制約のため、各運搬車は1つのバッテリーを持つことしかできない。本文は、バッテリー容量と積載量という二重の制約の下で、最適なサービス順と経路を追求し、参考のためにコストを計算する。

再配置のモデルを改善するために、次のような仮説を立てた：地域全体のシェアサイクルの総数が一定であり、地域全体のシェアサイクルが外部から回収されることなく、新たにシェアサイクルが利用されることもない、すべてのシェアサイクルが利用可能であると仮定する；各運搬車が同じモデルで、同じ積載量を持ち、同じバッテリー容量を持つと仮定する；各シェアサイクルが同じモデルで、同じ重量で、同じように使用できると仮定する；各運搬車の燃料消費能力はすべてのケースで同じであり、走行距離と積載量に線形かつ正の相関関係があると仮定し、都市部の複雑な道路や、運搬車自体の新旧程度による影響は考慮しない。再配置の過程における各運搬車は、完全に充電された状態で出発し、再配置の全過程において充電は禁止されていると仮定する。

3.2 モデルの構築

再配置経路モデルの全体的な枠組みは、従来のVPR問題のモデルをベースとしているが、シェアサイクルの再配置問題の特性や実用を考慮し、本文の再配置のモデルにおける制約条件には、各タスクにおける運搬車の最大走行距離、走行中の燃料消費量、最大積載量、バッテリー容量、最大配車回数、シェアサイクルの積載量が増えるごとに消費量が増

加することを考慮している。コスト要因は需要量と往復の実距離を考慮した。実際の運用では、事業者はユーザーの要求を瞬時に満たし、短時間で可能な限り再配置作業を完了させる必要があるため、需要を優先する貪欲アルゴリズム、距離を優先する貪欲アルゴリズム、および省エネを優先する貪欲アルゴリズムの 3 つの貪欲戦略を作成した。私のデータセットの更新はリアルタイムで行われ、つまり単位時間ごとにシェアサイクルが更新される。基本的な VPR 問題とは異なり、シェアサイクルの再配置問題では、各ポート点は需要側であると同時に、供給側でもなりうる。従って、我々の運搬車はシェアサイクルの積み下ろしを同時に行うことができる。特定の配車ステーションを持たないため、運搬車はタスクを完了した後、スタート地点に戻る必要がない。

3.2.1 シンボルの設置

具体的は以下のようになる。

表 1

| シンボル | 説明 |
|---|---|
| n | ここではポート点の数は 150 |
| $d_{ij}$ | ポート点 i と j の間の距離 |
| $c_i$ | ポート点 i の最大なシェアサイクルの数。 |
| $b_i$ | ポート点 i の現在のシェアサイクルの数。 |
| $B_i$ | ポート点 i のシェアサイクルのベ |

|  | ースライン。 |
|---|---|
| V | ここでは運搬車の最大数は 20。 |
| C | ここでは運搬車 1 台あたりの最大積載量は C =トラックの積載量/シェアサイクルの重量= 1150/30.2=38。 |
| D | ここでは運搬車 1 台あたりの最大走行距離は 270000 メートル。 |
| E | ここでは各運搬車のバッテリー容量は 38,700 Kwh |
| e | ここでは各運搬車が 1 メートル進むたびに消費する電力は 0.1433 Kwh |
| e′ | 運搬車 1 台がシェアサイクルを 1 台積むごとに、1 メートル歩くごとに 0.00327 Kwh の電力を消費する |
| P | ここではワット時あたりの電気料金は 0.03356 円 |
| F | ここでは運搬車 1 台あたりの固定費は、2,740 円。 |
| f | 運搬車がシェアサイクルを 1 台搭 |

| | 載して出発するあるいはポート点から1台回収すると、1単位あたり100円の追加コストがかかる。 |
|---|---|
| $x_{ijk}$ | 運搬車kがデータ点iからデータ点jまでの場合は1であり、そうでない場合は0。 |
| $q_{ik}$ | 運搬車kのポート点iへの積載量。 |

3.2.2 データセットの更新

1). ポアソン分布：確率質量関数（PMF）は次のようになる。

$$P(X = k) = \frac{\lambda^k e^{-\lambda}}{k!}$$

ここで、Xはある時間または空間における事象の発生回数を表す確率変数、kは特定の回数、λは1単位時間におけるランダム事象の平均発生率、eは自然対数の底である。このモデルでは、ポアソン分布を用いて、各時間帯における各駅の現在の自転車数をシミュレートしている。パラメータλは、時間帯とサイトの総容量に基づいて計算され、その時間帯における各サイトの自転車の平均利用率を表す。そして、ポアソン分布を使って生成された乱数が、現在の自転車数を表している。

2). 多項分布：多項分布の確率質量関数（PMF）は次のとおりです。

$$P(X_1 = n_1, X_2 = n_2, ..., X_k = n_k) = \frac{n!}{n_1! n_2! ... n_k!} p_1^{n_1} p_2^{n_2} ... p_k^{n_k}$$

このうち、$X_i$は i 番目のイベントが発生する確率を表す確率変数、$n_i$は i 番目のイベントの特定の発生回数、$p_1$ は i 番目のイベントの確率を表します発生；n はすべてのイベントの発生の合計数です。

このモデルでは、多項分布を使用してステーションをランダムに選択し、このステーションの過剰な自転車数を他のステーションに分配します。ここでの多項分布のパラメータは、全サイトの総容量に対する各サイトの総容量の比であり、各サイトが選択される確率を示す$c_i / \Sigma c_i$。である。次に、多項分布を使用して生成された乱数を使用して、自転車の数を増やすために選択されたステーションを表します。$B_i = c_i *$ percent，$B_i$はステーションの総容量，percent は自転車利用率です。地铁站附近的点和企业学校附近的点会赋予另外的值并进行泊松分布，而地下鉄駅付近の点と企業学校付近の点は加算値を付与してポアソン分布で分布し、地下鉄駅付近と企業学校付近の点は除外して直接計算します。事前に設定された percent。各ステーションの余剰自転車台数は、ステーションの総収容台数と現在の自転車台数に基づいて計算されます，式：diff = $c_i$ − $b_i$。

最終的に計算されたデータセットを下の図に示します。

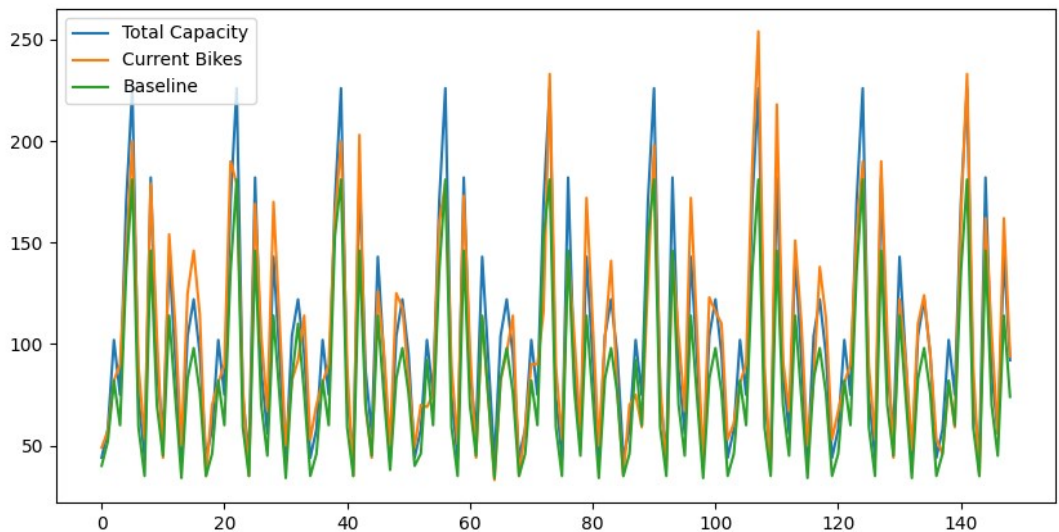

図6 ある瞬間におけるシェア自転車台数の推移

参考として公開データセットを使用しており、青い線は各駐車場の総収容台数、オレンジ色の線は現在のシェア自転車台数、緑の線はベースライン値を表しています。

単位時間ごとに調整しており、現状の自転車台数は全体の79%、基準値は全体の80%となっております。 特別ポイントの優遇を行っております 具体的には、地下鉄の駅近くに11ポイント、学校や企業の近くに23ポイントを設けており、自転車の使用量が多いなどの優遇をさせていただいております。自転車総数の68%、ベースラインは自転車総数の90%であり、これらの比率は時点によって異なります。 各拠点で一定の割合で減らされたシェア自転車の台数は、多項分布に従って他の拠点にランダムに配分され、最終的に各拠点のシェア自転車台数は総収容台数の1.2倍を超えることはできません。

3.2.3 モデリング

各データ ポイントは 1 つのディスパッチャによってのみ処理できます。

$$\sum_{k=1}^{v}\sum_{j=1,j\neq i}^{n} x_{ijk} = 1, \forall i \in \{1,2,\dots,n\}$$

配車車両の積載量は最大積載量を超えることはできません。

$$\sum_{i=1}^{n} q_{ik} \leq C, \forall k \in \{1,2,\dots,v\}$$

配車車両の走行距離は、以下の最大走行距離を超えることはできません。

$$\sum_{i=1}^{n}\sum_{j=1}^{n} d_{ij}x_{ijk} \leq D, \forall k \in \{1,2,\dots,v\}$$

配車車両の消費電力はバッテリー容量を超えることはできません。

$$\sum_{i=1}^{n}\sum_{j=1}^{n} (e + e'q_{ik})d_{ij}x_{ijk} \leq E, \forall k \in \{1,2,\dots,v\}$$

データ ポイント内の共有自転車の数は、次の要件を満たす必要があります。

$$b_i - \sum_{k=1}^{v} q_{ik} > B_i, \forall i \in \{1,2,\dots,n\}$$

　本論文では、ルートの長さが電気配車車両の消費電力に与える影響だけでなく、負荷の違いが燃料配車車両の消費電力に与える影響も調査し、総合的にルートを計画する。　現時点では、電気トラックの電力消費に影響を与える要因が多く、各要因が車両の電力消費に与える影響を区別することが難しいため、ルート最適化の目標として車両の消費電力に関する研究はほとんどありません。この種の研究の難しさの 1 つは、経路、負荷、消費電力の関係を計算することにあります。　現在、一部の学者が車両問題について研究を行い、さまざまなモデルを

確立し、さまざまな計算方法を提案しています。 本稿では経路計画を主に研究するため、既存手法の中でも比較的単純な「負荷推定手法」を採用する。 この方法は、車両の走行距離と負荷が消費電力を決定すると考えており、負荷と消費電力の間には線形関係があり、その関係は次のとおりです。

3.2.4 貪欲な戦略

最短距離を優先する、より大きな需要を満たすことを優先し、需要が同じであれば最短距離で満足させる、という3つの欲深い戦略があります。消費電力の最小化を優先します。

配車担当者が配送ステーションからスタートし、まずシェア自転車が多すぎる場所でシェア自転車を回収し、その後シェア自転車が不足している場所に配送するという最短経路を優先した欲張りな戦略。車両に積まれている自転車がすべて降ろされているか、降ろされていない場合は、さらに自転車を回収する必要があるため、車両は配送ステーションに戻ります。 具体的な手順は次のとおりです。

1.固定費、積載費、最大走行距離、配車車両の最大積載量などの定数を設定します。

2. 現在のシェア自転車台数とベンチマーク番号を含むシェア自転車ステーションデータをロードし、ステーション間の距離である距離マトリクスデータをロードします。

3. シェア自転車が多すぎるステーション（現在のシェア自転車数が基準数より多い）とシェア自転車が不足しているステーション（現

在のシェア自転車数が基準数より少ない）を見つけ、これらを保管するための 2 つのキューを作成します。駅を選択し、距離が小さいものから大きいものまで並べ替えます。

4.配車車両を初期化し、最初は1台のみ配車します。 ループに入ります。処理する必要があるサイトがまだある場合は、操作を続行します。

4.1 ステーションが多すぎる場合は、キューから 1 つを取り出し、配車車両の走行ルートに追加し、車両の積載量、積載量変化、積載率、走行距離、エネルギー消費量、その他の属性を更新します。

4.2 配車車両にシェア自転車が積載されており、シェア自転車が不足しているステーションがある場合、不足しているキューから1台を取り出して配車車両の走行ルートに追加し、積載量、積載変更、積載率を更新する車両の、走行距離、エネルギー消費量、その他の属性。

4.3. 車両が満載の場合、または移動する共有自転車がなくなった場合、車両は配達ステーションに戻ります。

5. 処理が必要なシェア自転車があり、配車台数が上限に達していない場合は、新たな配車車を配車します。

6. 各配車車両の走行ルート、積載量、積載変化、積載率、走行距離、エネルギー消費量、コスト、および配車車両全体の総距離、総エネルギー消費量、総コストを計算し出力します。

このプロセスでは、配車担当者は毎回最寄りのステーションのみを選択し、車両の数が最大に達するか需要が完全に満たされるまで、毎

回新しい車両のみを配車するという貪欲な戦略を採用しています。

消費電力が最小のパスを使用する場合も手順は同様です。経路長を消費電力に変換します。より大きな需要に基づいた貪欲な戦略が採用される場合、プログラムの 3 番目のステップ「キューを作成した後に距離で並べ替える」を「最大の需要を持つデータ ポイントが最初にランク付けされ、需要が同じ場合は、ユーザーのニーズを満たすことを優先するという目的を達成するため。なぜなら、共用駐輪場に不足する自転車が増えれば増えるほど、利用者の需要は高まるからです。

したがって、このプロセスでは、需要が多く、距離が短い駅を優先してできるだけ早く需要に応え、走行距離を短縮し、車両数が上限に達するまで一度に新しい車両を送り込むという欲張りな戦略をとります。最大値に達するか、需要が完全に満たされます。

4 論文の結果と分析

4.1 パラメータの設計

このペーパーでは、VRP 問題を解決するための異なる戦略を持つ 3 つの貪欲アルゴリズムを紹介し、これら 3 つのアルゴリズムを使用してスケジューリング パス最適化モデルを解決し、それらを比較します。最適経路を解く前に、共用駐輪場と共用駐輪場間の距離行列をメートル単位で整理する必要がある。

以前に、駐車ポイントの数が 150 個の場合の各ポイント間の距離行列を示しました。最初のポイントは開始点を表し、ノード 0 〜 150

は 151 個の駐車ポイントの位置を表します。 各駐車ポイントからの距離が検索に影響を与えるのを避けるため、これらの距離は 0 に設定されます。

4.2 実験結果の解析

貪欲アルゴリズムとは、問題を解決する際の目標として局所的な最適解を追求することを指します。 私たちは 3 つの貪欲な戦略を設計しました。1 つは、現在のノードまでの最短距離に基づいて次に訪問するノードを見つけることです。 もう 1 つは、共用駐輪場の需要に基づいており、需要が安定している場合は距離係数を参照します。 3 番目は、ノード間の最短の電力消費に基づいています。 スケジューリングモデルでは、目的関数は全体コストの最小化を追求しており、その中で物流トラックの性能、共用駐輪場間の往復距離、輸送コスト、共用自転車の積み降ろしコストなどは実調査に基づいて決定されています。 シェア自転車のスケジューリングルートは、最適化プロセスを介さずにノードの検索時に直接形成されるため、優れたルートを探す際には、シェア自転車駐車場の需要にさらに注意を払う必要があります。 すべての配車トラックは、積載量とバッテリー寿命の点で妥当な制限内にあります。

最短距離に基づくグリーディアルゴリズムの総距離は 209575.7 メートル、スケジューリングコストは 66293.3 円、消費電力は 121431.8Wh となる。 オペレーターは、式 [0、2、3、0、5、0、10、16、0、23、29、37、0、49、59、0、71、91、0、 104、119、0、139、

12、70、106]、[0、4、6、11、0、17、24、30、0、38、50、0、60、73、0、92、105、0、122、140、0]、[0、9、13、18、25、0、31、0、39、0、51、61、74、0、93、107、123、0、141、14、72、114]、[0、15、19、0、26、32、0、40、53、0、62、0、78、94、0、108、124、142、21、76、121]、[0、20、27、0、33、0、41、54、0、63、80、0、95、109、126、0、143、22、77、0]、[0、28、34、0、42、0、55、0、64、81、96、0、110、127、144、36、79、125]、[0、35、44、56、65、0、82、98、111、0、128、145、0]、[0、48、57、67、83、0、99、0、112、129、0、146、43、84、130]、[0、58、68、86、0、100、113、0、132、147、0]、[0、69、88、101、0、115、133、0、148、45、87、131]、[0、90、0、102、116、0、135、0]、[0、103、117、0、136、0、46、89、0]、[0、118、0、137、1、47、0]、[0、138、8、66、97、134]を使用して共有自転車ネットワークにアクセスします。

最大デマンドに基づくグリーディアルゴリズムの総距離は 233338.7 メートル、スケジューリングコストは 67690.8 円、消費電力は 1397319.1Wh となります。 オペレーターは、[0, 39, 0, 42, 0, 55, 0, 90, 0, 11, 0, 99, 0, 28, 138, 0, 32, 74, 0, 103 , 88、0、143、35、93、0]、[0、25、0、123、0、73、0、113、0、116、0、50、19、0、60、91、0、10、137、0、15、58、26、1、106、46、125]、[0、5、0、62、0、59、0、133、0、49、94、0、104、105、0、27、17、0、81、92、107、47、134、70、130]、[0、147、0、140、0、33、

0、63、111、0、3、122、0、61、24、0、82、30、127、87、12、79、131]、[0、31、0、118、0、80、128、0、20、139、0、78、41、0、98、110、13、0、14、84]、[0、135、0、148、145、0、40、2、0、95、68、0、115、6、9、8、36、97]、[ 0、65、67、0、53、83、0、112、69、0、132、64、38、21、43、114]、[0、96、54、0、141、146、0、100、23、124、0、44、22、0]、[0、71、16、0、117、102、0、108、101、48、66、45、121]、[0、86、129、0、119、109、142、0、57、76、0]、[0、37、126、144、0、4、89、0]、[0、136、51、29、0]、[0、18、[34、72、0]、[0、56、77、0]のパスは共有自転車ネットワークにアクセスします。

最短消費電力に基づくグリーディアルゴリズムの総距離は 196212.9 メートル、スケジューリングコストは 62724.5 円、消費電力は 106471.3Wh となる。オペレーターは、式 [0, 1, 14, 0, 4, 0, 6, 55, 0, 15, 104, 0, 24, 48, 0, 37, 60, 76, 0, 87、138、135、0]、[0、2、84、0、7、17、57、0、20、112、108、0、39、70、0、73、0、91、111、0、123、126、32、74、0]、[0、5、0、8、0、11、71、0、21、0、25、0、31、0、41、56、0、75、103、145、0、129、125、134、100]、[0、9、63、0、18、34、33、0、35、89、0、50、77、105、113、0、130、52、58 , 98]、[0, 12, 53, 0, 22, 30, 106, 0, 42, 0, 54, 110, 0, 93, 0, 115, 139, 114, 146, 43]、[0、19、67、0、27、40、79、0、59、0、78、109、147、0、133、132、0]、[0、23、29、38、0、44、122、148、140、0、116、0、136、3、82、68、

0]、[0、28、45、128、0、61、102、97、127、0、141、0]、[0、36、90、0、62、101、0、94、144、0、143、47、72、80]、[0、49、0、65、107、0、95、137、51、46、16、92]，[0, 66, 0, 81, 96, 142, 120, 13, 10, 0], [0, 86, 121, 0, 117, 0]，[0, 99, 131, 0]，[0、118、0]を使用して共有自転車ネットワークにアクセスします。

表2 貪欲なアルゴリズム戦略の比較

|  | 最短距離に基づく貪欲アルゴリズム | 最大需要に基づく貪欲アルゴリズム | 最短の消費電力に基づいた貪欲なアルゴリズム |
| --- | --- | --- | --- |
| 最短距離 | 209575.7 メートル | 233338.7 メートル | 196212.9 メートル |
| 消費電力 | 287365.6 ワット時 | 296855.5 ワット時 | 273257.1 ワット時 |
| 最低コスト | 66293.3 円 | 67690.8 円 | 62724.5 円 |
| 配車台数 | 14 台 | 14 台 | 14 台 |
| 消費電力 | 121431.8 Wh | 1397319.1 Wh | 106471.3 Wh |

5．結論は

この表から、それぞれの貪欲アルゴリズムが異なるメトリックにおけるパフォーマンスを評価することができます。以下に、各アルゴリズムの特性とその結果に対する解釈を提供します。

最短距離に基づく貪欲アルゴリズム：このアルゴリズムは、最短の距離を優先して選択することにより、全体の移動距離を最小化します。しかし、このアルゴリズムはエネルギー消費や総コストには最適化されていません。これは、一部の近いルートが高いエネルギー消費を伴う可能性があるためです。

最大需要に基づく貪欲アルゴリズム：このアルゴリズムは、最も需要の高い地点を優先します。これにより、可能な限り多くの自転車を効率的に移動することができます。しかし、このアルゴリズムもエネルギー消費や総コストには最適化されていません。最大の需要を持つ地点が必ずしも最もエネルギー効率的なルートであるとは限らないからです。

最短の消費電力に基づいた貪欲なアルゴリズム： このアルゴリズムは、エネルギー消費を最小化するルートを優先します。その結果、最小の消費電力と最低コストを達成することができます。しかし、このアルゴリズムは距離には最適化されていません。一部のエネルギー効率の良いルートは、より長い距離を必要とするかもしれません。

以上の結果から、最適なアルゴリズムは目的によります。エネルギー消費とコストを最小化したい場合、消費電力に基づく貪欲アルゴリズムが最良の選択となります。しかし、移動距離を最小限に抑えたい場合は、

最短距離に基づく貪欲アルゴリズムを選択することが最善かもしれません。

参考文献.

**ドコモとドコモ・バイクシェアの使用車種**
https://www.cycling-ex.com/2015/09/minatoku_share.html
https://www.yodobashi.com/product/100000001002947805/

ドコモとドコモ・バイクシェアにおける最適配置作業用運搬車

https://www.nextmobility.jp/economy_society/demonstration-of-linking-follow-fly-ev-trucks-and-bicycle-sharing20230130/

運搬車使用車種

https://folofly.com/product/

電気料金 https://www.globalpetrolprices.com/Japan/electricity_prices/

駅 https://www.tokyometro.jp/station/index.html

港区自転車シェアリングポート

https://catalog.data.metro.tokyo.lg.jp/dataset/t131032d0000000011/resource/6c58af42-2b29-4e18-878d-08a64d69a590

前田 謙太郎,"バイクシェアリングにおける自転車再配置車両の経路決定手法に関する研究,"日本オペレーションズ・リサーチ学会 学生論文賞受賞論文,2017.11

LIU Xin-yu, CHEN Qun, "An Optimization Model for Bike Repositioning in Bike-sharing Systems Considering Both Demands for Borrowing or Returning Bikes and Costs of Repositioning Operations," China J. Highw. Transp. 10.19721/j.cnki.1001-7372.2019.07.016